\documentclass[russian]{article}
\usepackage{textcomp}
\usepackage[ukrainian,english]{babel}
\usepackage[T2A]{fontenc}
\usepackage{indentfirst}

\usepackage[utf8x]{inputenc} 

\usepackage{amssymb}

\usepackage[utf8x]{inputenc} 
\usepackage{amsmath}
\usepackage{color}
\usepackage{graphics}
\usepackage{epsfig}
\usepackage[all,knot]{xy}
\xyoption{arc}

\def\mapr#1{\smash{\mathop{\buildrel{#1}\over\longrightarrow}}}

\newtheorem{theorem}{Theorem}

\def\qed{\hfill\vrule width2mm height2mm depth2mm}

\def\proof{{\bf Proof.}}

\def\b{{\bf b}}

\def\e{{\bf e}}

\def\0{{\bf 0}}
\def\1{{\bf 1}}

\def\b#1{{\bf #1}}

\def\cD{{\cal D}}

\def\cZ{{\cal Z}}

\def\ad{{\hbox{\bf ad}}}

\def\Der{\mbox{\bf Der }}

\def\im{{\hbox{\bf Im}\;}}

\def\Out{{\hbox{\bf Out}\;}}

\def\fcolor{\color{black}}

\def\beq#1{\fcolor\begin{equation}\label{#1}}
\def\eeq{\end{equation}\color{black}}

\def\b${\color{blue} $}
\def\e${$ \color{black}\hskip-0.05cm}

\date{}

\title{Hochshild's method for describing the Mackenzie obstruction for
construction of transitive Lie algebroid.}

\author{Gasimov V.A.}

\begin{document}
\maketitle

{\textbf{Introduction.}}

Extensive development of noncommutative geometry began in the middle of the last century. The origins of noncommutative geometry go back to the Gelfand-Naimarck theorem on a one-to-one correspondence between the category of compact topological spaces and the category of commutative C*-algebras. The question of which geometric images can correspond to noncommutative C*-algebras has led to the creation of noncommutative geometry. Noncommutative geometry studies the constructions related to noncommutative algebras which, for the commutative case, have a well-defined geometric interpretation. The idea of noncommutative geometry has enabled joining the variety of topology representations and methods, differential geometry, functional analysis, representation theory, etc. In particular, groups and Lie algebras, as well as groupoids and Lie algebroids, serve as examples of noncommutative structures. 

This work consists of an introduction and four parts. In the first  part provides the necessary definitions and facts related to this work. The theory of groupoids and Lie algebroids is in many ways analogous to the theory of groups and Lie algebras. There is an analogy between the connection of Lie algebroids with Lie groupoids and the connection of Lie algebras with Lie groups. Similarly, any Lie group can be associated with a Lie algebra and any Lie groupoid can be associated with a Lie algebroid \cite{Mck-2005}. The structure of Lie algebroid reflects various aspects of smooth manifolds associated with differential operators and dynamical systems on manifolds. Smooth manifolds and various differentiation operations are very closely related to Lie algebroids. The study of transitive Lie algebroids was started in the book “General Theory of Lie Groupoids and Lie Algebroids” by Professor Mackenzie.

Some results obtained in the Mackenzie study of the obstruction for constructing a transitive Lie algebroid were presented by the author together with Professor A. Mishchenko and Chinese mathematician Xiaoyu Li at the International Mathematical Congress in Seoul   in 2014 \cite{qv-2019}. The same result was reported at the Topological Conference in Navpaktos (Greece) \cite{qv1-2019}, \cite{M.A.1}.

 The second part examines extensions of Lie algebras.
 The exact sequence of a transitive Lie algebroid can be considered as an extension of a Lie algebra with the help of another Lie algebra. More precisely, the algebra of vector fields $\Gamma ^{\infty }\left( TM \right)
\ $  is extended to Lie algebra $\Gamma ^{\infty }\left( A\right)
\ $
  with the help of Lie algebra $\Gamma ^{\infty }\left( L\right)
  \ $ , which is called the extension nucleus \cite{hochschild-1954-1}. 
The difference between the extensions of Lie algebras and the exact Atya sequence of a transitive algebroid is that the structure of the commutator bracket $\left\{ \cdot ,\cdot \right\} $  is subordinated to the Leibnitz additional property-identity with respect to multiplication by smooth functions. A homotopy classification of transitive Lie algebroids is carried out in this paper \cite{M.A.2}. This problem is solved on the basis of a category approach and boils down to verifying the triviality of Mackenzie’s obstruction for the existence of Lie algebroid.

Earlier in the paper \cite{qv-2019}, the triviality of Mackenzie's obstructions for the Lie algebras  ${\mathfrak g}/Z_{\mathfrak g}$  was shown. The question of calculating the obstruction for coupling was considered for the simplest cases of finite-dimensional Lie algebras. The paper shows \cite{qv-2019} that if the Mackenzie obstruction is trivial for the layering of algebras  ${\mathfrak g}_{1}$ and ${\mathfrak g}_{2}$ , then the Mackenzie obstruction will also be trivial for a direct sum of algebras ${\mathfrak g}_{1} \bigoplus {\mathfrak g}_{2}$ . To solve the problem of classification for transitive Lie algebroids with a fixed adjoint finite-dimensional Lie algebra, it is also necessary to study the automorphism group ${\bf Aut(\mathfrak g)}^{\delta}$ along with the new topology. At the conference in Krakow (2016) \cite{qv2-2019}, the author made a report on the topic of automorphisms of generalized Heisenberg algebra. The paper \cite{qV.4} describes the automorphism group of the generalized Heisenberg algebra.

  The third part shows that according to Mackenzie's terminology \cite{Mck-2005}, homomorphism $\Xi$ is a coupling between the Lie algebra $L$ and $T$ only when both algebras are spaces of sections of some vector bundles. 

Thus, the technique developed in Hochschild's works \cite{hochschild-1954}, \cite{hochschild-1954-1} may be applied to the special case of transitive Lie algebroids.
In the fourth part of the paper \cite{hochschild-1954}, we will formulate the problem of describing the Mackenzie obstruction for constructing a transitive Lie algebroid. The Mackenzie obstruction for the given coupling is a three-dimensional cohomology class, which is represented by a closed differential form

\begin{equation*}
\left(d^{\nabla}\Omega\right) \left(X,Y,Z\right) \in\ \Gamma ^{\infty }\left( L\right). 
\end{equation*}

Therefore, the natural problem of computing and describing the cohomology class arises $%
Obs\left(L, \Xi \right) \in H^{3}\left( M;ZL;\nabla ^{Z}\right) $. In the case of Lie algebras in Hochschild's work \cite{hochschild-1954} it was shown that in such a set one can define the structure of a linear space, and the map defined by the Mackenzie obstruction is a linear monomorphism. Our problem stated as that to transfer the Hochschild construction to the case of transitive Lie algebroids and to prove a similar theorem. We show that the image

\begin{equation*} 
\im(UObs)\subset H^{3}(M;\cZ)
\end{equation*}

is a linear subspace. The results of this work were presented by the author at a conference dedicated to the memory of the outstanding mathematician Boris Sternin (2018) \cite {qv3-2018}.

\section{ Definition of transitive Lie algebroid.}

The Lie algebroid $A$ on a smooth manifold $M$ is a vector bundle $%
p:A\rightarrow M$ together with structure of the Lie algebra
$\left\{ \cdot ,\cdot \right\} $ on the space $\Gamma ^{\infty
}\left( A;M\right) $ of smooth sections and a mapping of bundles
$a:A\rightarrow TM$ that is called the anchor, such that

(i) the induced mapping $\Gamma(a):\Gamma ^{\infty }\left( A;M\right)
\rightarrow \Gamma ^{\infty }\left( TM;M\right) $ is the
homomorphism of  Lie algebras;

(ii) for any sections $\sigma, \tau \in\Gamma^{\infty }\left(
A;M\right) $ and smooth function $f\in C^{\infty }\left( M\right)
$ the Leibnitz identity is fulfilled with respect to the operation
of multiplication of the sctions by the function
\begin{equation*}
\left[ \sigma ,f\tau \right] =f\left[ \sigma ,\tau \right]
+a\left( \sigma \right) \left( f\right) \tau .
\end{equation*}

The Lie algebroid $A$ is said to be transitive if the anchor $a$ is
fiberwise surjective. For transitive Lie algebroid one has the exact Atiyah sequence
\begin{equation*}
0\rightarrow L\overset{j}{\rightarrow }A\overset{a}{\rightarrow }%
TM\rightarrow 0.
\end{equation*}

The Mackensie obstruction is a three dimensional class of
cohomologies \cite{Mck-2005} whose
triviality is provided by the existence and construction of the
transitive Lie algebroid on the manifold
 $M$ if we are given the set a data:

1) The local trivial bundle $L$ will typical fiber isomorphic to
the finite-dimensional Lie algebra $g$ and structural group of all
automorphisms $Aut\left( g\right) $ of algebra $g$, denoted by
$LAB$.

2) Coupling between the tangential bundle $TM$ and bundle $LAB$ in
the form of homomorphism
\begin{equation*}
\Gamma ^{\infty }\left( \Xi \right) :\Gamma ^{\infty }\left(
TM\right) \rightarrow \Gamma ^{\infty }\left(
\mathcal{D}_{out}\left( L\right) \right)
\end{equation*}
of the space of sections with respect of the infinite dimensional Lie algebra structures.

According to this set of data, in Mackenzie's book \cite{Mck-2005} three-dimensional class of cohomologies is
structured in the following way. The exact sequence of the Lie
algebroid is considered:
\begin{equation*}
0\rightarrow ZL\rightarrow L\overset{ad}{\rightarrow }\mathcal{D}%
_{der}\left( L\right) \overset{\natural }{\rightarrow }\mathcal{D}%
_{out}\left( L\right) \rightarrow 0.
\end{equation*}
The coupling $\Xi $ as a homomorphism of bundles rises to the connection $%
\nabla $ in the diagram \hspace{1mm}
$$
\xymatrix{
0\ar[r]& ZL\ar[r]& L\ar[r]^{ad}&\mathcal{D}_{der}(L)\ar[r]^{\natural}& \mathcal{D}_{out}(L)\ar[r]& 0\\
&&&TM\ar[u]^{\nabla}\ar[ur]_{\Xi}&&
}
$$

Since the coupling $\Xi$ generates the homomorphism of infinite-dimensional Lie
algebras, i.e. the condition%
\begin{equation*}
\Xi _{\left[ X_{1},X_{2}\right] }=\left[ \Xi _{X_{1}},\Xi
_{X_{2}}\right] ,
\end{equation*}
holds, then the curvature tensor of the coupling $\Xi$ is trivial
\begin{equation*}
R^{\Xi }\left( X_{1},X_{2}\right) =\left[ \Xi _{X_{1},}\Xi
_{X_{2}}\right] -\Xi _{\left[ X_{1},X_{2}\right] }\equiv 0
\end{equation*}%
i.e. for the curvature tensor $R^{\nabla }$ the following
condition holds:%
\begin{equation*}
\natural \circ R^{\nabla }=R^{\Xi }=0.
\end{equation*}

This means that the curvature tensor $R^{\nabla }$ as the map in the diagram
$$
\xymatrix{
0\ar[r]& ZL\ar[r]& L\ar[r]^{ad}&{\Der}(L)\ar[r]^{\theta}& {\Out}(L)\ar[r]& 0\\
&&&\Lambda^{2}TM\ar[lu]^{\Omega}\ar[u]^{R^{\nabla}}\ar[ur]_{R^{\Xi}=0}&&
}
$$
rises to the homomorphism $\Omega $, since the sequence %
\begin{equation*}
0\rightarrow ZL\rightarrow L\overset{ad}{\rightarrow
}\mathbf{Der}\left( L\right) \overset{\theta}{\rightarrow
}\mathbf{Out}\left( L\right) \rightarrow 0
\end{equation*}
is also exact.

The differential $d^{\nabla }\Omega $ of the form $\Omega $ is
included in the commutative diagram \hspace{1mm}
$$
\xymatrix{
0\ar[r]& ZL\ar[r]^{i_{ZL}}& L\ar[r]^{ad}&{\Der}(L)\ar[r]^{\theta}& {\Out}(L)\ar[r]& 0\\
&&\Lambda^{3}TM\ar[u]^{d^{\nabla}\Omega}\ar[ur]_{d^{\nabla^{der}}R^{\nabla}}&&&
}
$$
For the curvature tensor the so called Bianchi
differential identity  that may be expressed in the form of the
equality $d^{\nabla der}R^{\nabla }=0$ holds. Therefore,
the differential $d^{\nabla }\Omega $ takes
values in the bundle $ZL$:%

$$
\xymatrix{
0\ar[r]& ZL\ar[r]^{i_{ZL}}& L\ar[r]^{ad}&{\Der}(L)\ar[r]^{\theta}& {\Out}(L)\ar[r]& 0\\
&&\Lambda^{3}TM\ar[u]^{d^{\nabla}\Omega}\ar[ur]_{d^{\nabla^{der}}R^{\nabla}}
\ar[ul]^{i_{ZL}^{-1}(d^{\nabla}\Omega)}&&&
}
$$

As the connection $\nabla $ is invariant on the subbundle $ZL$ we have the connection $\nabla ^{Z}$
induced on the subbundle $ZL$. The connection $\nabla ^{Z}$ has the
trivial curvature tensor
\begin{equation*}
R^{\nabla^Z}\equiv 0.
\end{equation*}
Moreover
\begin{equation*}
d^{\nabla^Z}\left( d^{\nabla }\Omega \right) \equiv 0,
\end{equation*}
therefore the bundle $ZL$ is a flat bundle. The form
$i_{ZL}^{-1}(d^{\nabla }\Omega)$ is a closed form in the space of 3
forms $\Omega^{3}\left(M;ZL,\nabla^{Z}\right)$.

The cohomologies  class  $[i_{ZL}^{-1}\left( d^{\nabla }\Omega
\right)]$  in the group $H^{3}\left( M;ZL,\nabla
^{Z}\right) $, given by the form $i_{ZL}^{-1}\left( d^{\nabla }\Omega
\right) $ is denoted by

\begin{equation*}
Obs\left(L, \Xi \right)=[i^{-1}\left( d^{\nabla }\Omega
\right)] \in H^{3}\left( M;ZL;\nabla
^{Z}\right)
\end{equation*}
and is called the Mackenzie obstruction for the coupling $\Xi $.

According to theory 7.2.13 from Mackenzie is book \cite{Mck-2005}
the cohomologies class $Obs\left(L, \Xi\right) \in
H^{3}\left( M;ZL;\nabla ^{Z}\right) $ is dependent only on the
coupling $\Xi $ and is independent of the choice of the connection
$\nabla$ and of the form $\Omega$.

This obstruction $Obs(L, \Xi)\in H^{3}(M;ZL;\nabla^{Z})$ equals zero if and only if the pair $\left(
\nabla ,\Omega \right) $ may be chosen such that $d^{\nabla }\Omega
\equiv 0$, i.e. forms the structure of the transitive Lie
algebroid for the given coupling $\Xi $.

Thus, for classifying transitive Lie algebroids with fixed
associated finite-dimensional Lie algebra bundle $L$ it is required to
check for which couplings $\Xi$  between bundles $TM$ and $L$,
$\Xi:TM\mapr{}\cD_{out}(L)$,
the cohomologies class
$Obs(L,\Xi)\in H^{3}(M;ZL;\nabla^{Z})$
is trivial. Therefore, there arises a
natural problem of calculation and description of the cohomologies class

\begin{equation*}  
Obs\left(L, \Xi \right) \in H^{3}\left( M;ZL;\nabla ^{Z}\right). \end{equation*}

\section{Extension of Lie algebras.}

The exact Atiyah sequence of the transitive Lie algebroid may be
treated as the extension of a Lie algebra by means of another
Lie algebra. More exactly, algebra of vector fields $\Gamma
^{\infty }\left( TM\right) $ is extended to the Lie algebra
$\Gamma ^{\infty }\left( \mathcal{A}\right) $ by means of the Lie
algebra $\Gamma ^{\infty }\left( L\right) $ that is called the
extension kernel (see: \cite {hochschild-1954-1}).

The difference of extensions of Lie algebras from Atiya's exact
sequence of transitive algebroid is that the structure of
commutator brace $\left\{ \cdot ,\cdot \right\} $ is subjected to
additional property, to the Leibnitz identity relative to
multiplication by smooth functions on a manifold. This fact is
explained in Wallas's \cite{Walas-2007} doctor dissertation, where instead of
Lie transitive algebroids, the so called Lie-Rinehart algebras are
considered by passing from bundles to cuts spaces of these
bundles.

More exactly, let $A$ be an arbitrary commutative algebra (over
the ring of scalars $R$), let $T=\mathbf{Der}\left( A\right) $ be
the Lie algebra of derivations of algebra $A$ that is the modulus
over algebra $A$. The algebra of smooth functions
\begin{equation*}
A=C^{\infty }\left( M\right)
\end{equation*}
on the smooth
manifold $M$
is a standard model  for the commutative algebra $A$. Then $T$ is a space of vector fields that is
isomorphic to the space of sections of tangent bundle
\begin{equation*}
T=\Gamma ^{\infty }\left( TM\right) .
\end{equation*}

In the model case when the manifold $M$ is compact, the algebra $T$ is a finitely-generated projective
modulus over algebra $A$.

\section{Definitions and Hochschild's terminology.}

The $A$-module $E$ equipped with commutator brace $\left\{\cdot ,\cdot
\right\}$ that defines a structure of Lie algebra, and the
homomorphism of Lie algebras
\begin{equation*}
E\overset{a}{\rightarrow }T,
\end{equation*}
(called an anchor), by means of which the Leibnitz identity with respect to the operation of multiplication of the
element of $E$ modulus by the element of algebra $A$ holds, is called
the Lie-Rinehart algebra.

If the anchor $a$ is an epimorphism, then we get expansion of
Lie-Rinehart algebras in the form of Atiya's exact sequence.

\begin{equation*}
0\rightarrow L\rightarrow E\overset{a}{\rightarrow }T\rightarrow 0.
\end{equation*}

Since the module $T$ is projective, the Atiya's exact sequence is splitting and $L$ is projective $A$-module. The module $L$ ia called
the extension kernel
\begin{equation*}
L=\mathbf{Ker}\left( a\right) .
\end{equation*}

By Hochschild's terminology, the kernel $L$ is $T$-kernel as it
is equipped with additional structure of module over the algebra $T$.
But the structure is  not classic, in the form of adjoint action
of Lie algebra $A$ on the modulus $L$
by the following rule:%
\begin{equation*}
E\times L\overset{\mathbf{ad}}{\rightarrow }L:\left( e,u\right) \mapsto %
\left[ e,u\right] \in L,\ \ \ e\in E,\ \ \ u\in L.
\end{equation*}

In fact we can understand this action of algebra $E$ as
multivalued action of
the quotient algebra $T$:%
\begin{equation*}
T\times L\rightarrow L:\left( X,u\right) \mapsto \left[
a^{-1}\left( x\right) ,u\right] \subset L,\ \ \ X\in T,\ \ \ u\in
L.
\end{equation*}%

The action $\mathbf{ad}$ may be defined as a homomorphism of the algebra
$E$ into algebra $\mathbf{Der}\left( L\right) $ of  derivations of
the algebra $L$:
\begin{equation*}
\mathbf{ad}:E\rightarrow \mathbf{Der}\left( L\right) ,
\end{equation*}%
that in composition with projection onto external derivations $\Out(L)$ gives the
homomorphism%
\begin{equation*}
E\overset{\mathbf{ad}}{\rightarrow }\mathbf{Der}\left( L\right) \overset{%
\natural }{\rightarrow }\mathbf{Out}\left( L\right) ,
\end{equation*}%
that is trivial on the kernel $L$, i.e. gives the well-formed homomorphism%
\begin{equation*}
T\overset{\Xi }{\rightarrow }\mathbf{Out}\left( L\right) .
\end{equation*}%
Thus, additional structure of $T$-kernel on $L$ modulus is given
as homomorphism $\Xi $ of $T$ algebra into algebra
$\mathbf{Out}\left( L\right) $ of external derivations of algebra
$L$.

By the Mackenzie terminology (\cite{Mck-2005}) the homomorphism $\Xi$ is called a coupling between the Lie algebra $L$ and algebra $T$
but only in the case when both algebras are spaces od sections of some
vector bundles. Thus, to the special case of transitive Lie
algebroids one can apply the technique worked out in Hochschild's
paper (\cite{hochschild-1954}).

\section{Hochshild description of the obstruction set (for transitive Lie algebroids)}

Following the paper \cite{hochschild-1954}, we formulate the problem of description
of Mackenzie obstruction to construction of transitive Lie
algebroid by the following  set of data:
\begin{enumerate}
\item The $LAB$ bundle  $L$ on a smooth manifold $M$,
\item The coupling $\Xi$ between the bundle $L$ and tangent bundle $TM$.
\end{enumerate}

Thus, we consider a vector bundle $\left( LAB\right) L\rightarrow
M$ whose typical fiber is the finite-dimensional Lie algebra $g$,
while the structural group is the group $\mathbf{Aut}\left(
g\right)$ of automorphisms of the fiber $g$. Fix the coupling
$\Xi $ between the bundle $L$ and tangent bundle $TM$:%
\begin{equation*}
\Xi :TM\rightarrow \mathcal{D}_{out}\left( L\right) ,
\end{equation*}%
which is the homomorphism of correspondent algebras of sections.  The lifting of the
coupling $\Xi$ to a connection $\nabla $,%
$$
\xymatrix{
&L\ar[d]^{\ad}\\
&\cD_{der}(L)\ar[d]^{\natural}\\
TM\ar[r]^{\Xi}\ar[ru]^{\nabla}&\cD_{out}(L)\ar[d]\\
&0
}
$$
 generates a lifting of the curvature tensor $R^{\nabla}$ to a two-dimensional form
 $$
\xymatrix{
&L\ar[d]^{\ad}\\
&\Der(L)\ar[d]^{\theta}\\
\Lambda^{2}(TM)\ar[r]^{R^{\Xi}=0}\ar[ru]^{R^{\nabla}}\ar@/^1pc/[ruu]^{\Omega}&\Out(L)\ar[d]\\
&0
}
$$
that lifts the differential $d^{\nabla }\Omega $ to a three-dimensional form $%
U=i^{-1}\left( d^{\nabla }\Omega \right) $:%
$$
\xymatrix{
&&&ZL\ar[d]^{i_{ZL}}\\
&&&L\ar[d]^{\ad}\\
&&&\Der(L)\ar[d]^{\theta}\\
\Lambda^{3}(TM)\ar[rrr]^{d^{\Xi^{out}}R^{\Xi}=0}\ar[rrru]^{d^{\nabla^{der}}(R^{\nabla})=0}
\ar@/^1pc/[rrruu]^{d^{\nabla}\Omega}\ar@/^1.5pc/[rrruuu]^{U=i_{ZL}^{-1}(d^{\nabla}\Omega)}&&&\Out(L)\ar[d]\\
&&&0
}
$$

So, the differential form $U=i_{ZL}^{-1}\left( d^{\nabla }\Omega \right)
$ defines the cohomology class

\begin{equation*}
Obs\left( \Xi \right)=[U] =\left[ i_{ZL}^{-1}\left( d^{\nabla }\Omega
\right) \right] \in H^{3}\left( M;ZL\right)
\end{equation*}
with values in the flat bundle $ZL$ wherein the flat structure
is set up by restriction of the connection $\nabla$ to the connection $\nabla^{ZL}$, that covers
the coupling $\Xi $, in the subbundle $ZL$.

According to Hochschild (\cite{hochschild-1954}) we consider the set
$$Coup=Coup(\cZ)=\{(L,\Xi, \varphi),\}$$
where $\Xi:TM\mapr{}\cD_{out}(L)$ is coupling, $\varphi:ZL\mapr{}\cZ$ is an isomorphism which preserves  linear connections $\nabla_{ZL}\approx \nabla_{\cZ}$,
of all bundles
(LAB) equipped with the coupling $\Xi $ with the tangent bundle
$TM$, at that the module $\Gamma ^{\infty }(ZL)$ is fixed over the Lie algebra $\Gamma ^{\infty }\left( TM\right)$ of vector fields. In the case
of Lie algebras in Hochschild's work (\cite{hochschild-1954}, pp. 698) it was shown
that in such a set it is possible to set up the structure of
linear space, and the mapping determined by the Mackenzie
obstruction is a linear monomorphism.

So, we have a map
$$
UObs:Coup(\cZ)\mapr{}H^{3}(M;\cZ),
$$
that is defined by the composition
$$
UObs:Coup(\cZ)\mapr{Obs}H^{3}(M;ZL)\mapr{\varphi^*}H^{3}(M;\cZ)
$$

Our task is to take the Hochschild construction to the case of
transitive Lie algebroids and to prove the similar theorem: the
set of all Mackenzie obstructions for bundles $(LAB)$ $L$ with
coupling $\Xi $ with  the same center $ZL$ as a modules
over the algebra of vector fields $\Gamma
^{\infty }\left( TM\right) $ generates a linear subspace in the group $%
H^{3}\left( M;ZL\right) $.

\begin{theorem}
The  image $\im(UObs)\subset H^{3}(M;\cZ)$ is a liner subspace.
\end{theorem}
\proof

Given two elements $c_i=(L_i,\Xi_i,\varphi_i)\in Coup$, $i=1,2$ we construct we construct the third element $c_3=(L_3,\Xi_3,\varphi_3)\in Coup$ by the following:

Consider the direct sum of couplings $L_1\oplus L_2$,
$$\Xi_1\oplus\Xi_2:\cD_{out}(L_1\oplus L_2)\approx
\cD_{out}(L_1)\oplus\cD_{out}(L_2),$$ $$\varphi_1\oplus\varphi_2:ZL_{1}\oplus ZL_{2}\mapr{}\cZ\oplus\cZ.$$

Consider the liner bundle homomorphism
$$
\theta:\cZ\oplus\cZ\mapr{}\cZ
$$
that is defined by the formula:
$$
\theta(u,v)=\alpha u+\beta v, \quad (u,v)\in\cZ\oplus\cZ, \quad \alpha,\beta\in \mathbb{R}.
$$

Let $\cZ_0=\ker\theta\mapr{i_{\theta}} \cZ\oplus\cZ$ and

$$
\psi: \cZ_{0}\mapr{}L_{1}\oplus L_{2}
$$
be the composition
$$
\psi=i_{Z}\oplus i_{Z}\circ(\varphi_1\oplus\varphi_2)^{-1}\circ i_{\theta}:
\cZ_{0}\mapr{}\cZ\oplus\cZ\mapr{}ZL_{1}\oplus ZL_{2}\mapr{}L_{1}\oplus L_{2}
$$
Then put
$$
L_{3}=(L_{1}\oplus L_{2})/\psi(\cZ_{0}).
$$

It is clear that coupling $\Xi_{1}\oplus\Xi_{2}$ induces the coupling $\Xi_{3}$ such that
$$
UObs(L_{1}\oplus L_{2},\Xi_{1}\oplus\Xi_{2}, \varphi_{1}\oplus\varphi_{2})=
UObs(L_{1},\Xi_{1}, \varphi_{1})\oplus
UObs(L_{2},\Xi_{2},\varphi_{2})\in H^{3}(M; \cZ\oplus\cZ),
$$
and 
$$
\begin{array}{l}
UObs(L_{3},\Xi_{3}, \varphi_{3})=\theta \left(
UObs(L_{1},\Xi_{1}, \varphi_{1}),UObs(L_{2},\Xi_{2}, \varphi_{2})
\right)=\\\\=
\alpha UObs(L_{1},\Xi_{1}, \varphi_{1})+\beta UObs(L_{2},\Xi_{2}, \varphi_{2})
\in H^{3}(M; \cZ)
\end{array}
$$
\qed

\end{document}